\documentclass{article}
\usepackage{amsthm}
\usepackage{amsmath}
\usepackage{amssymb}

\title{Rigidity of certain solvable actions\\ on the sphere}
\author{Masayuki ASAOKA
\footnote{Partially supported by JSPS Grant-in-Aid for Young
Scientists (A).}\\
Department of Mathematics, Kyoto University
}

\def\RR{\mathbb{R}}
\def\ZZ{\mathbb{Z}}

\def\cU{{\mathcal U}}
\def\cV{{\mathcal V}}

\def\cM{{\mathcal M}}

\def\cD{{\mathcal D}}
\def\cS{{\mathcal S}}

\def\GL{{\rm GL}}

\def\ra{{\rightarrow}}
\def\st{{\;|\;}}

\def\del{\partial}

\def\sm{{\setminus}}

\def\ra{{\rightarrow}}
\def\st{{\;|\;}}

\def\hsp{{\hspace{3mm}}}

\DeclareMathOperator{\Hom}{Hom}

\DeclareMathOperator{\Diff}{Diff}

\DeclareMathOperator{\Img}{Im}
\DeclareMathOperator{\Ker}{Ker}

\newcommand\Pair[2]{{\langle {#1},{#2}\rangle}}

\theoremstyle{plain}
\newtheorem{thm}{Theorem}[section]
\newtheorem{prop}[thm]{Proposition}
\newtheorem{lemma}[thm]{Lemma}
\newtheorem{cor}[thm]{Corollary}

\newtheorem*{Mthm}{Main Theorem}
\newtheorem*{propV}{Proposition \ref{prop:vanishing}}

\theoremstyle{definition}

\theoremstyle{remark}

\begin{document}
\maketitle

\begin{abstract}
An analog of Baumslag-Solitar's group $BS(1,k)$ naturally acts
 on the sphere by conformal transformations.
The action is not locally rigid in higher dimension,
 but exhibits a weak form of local rigidity.
More precisely, any perturbation preserves a smooth conformal structure.
\end{abstract}

\section{Introduction}
Over the last two decades,
 it has been found that
 many smooth actions of discrete groups exhibit local rigidity.
Most of known examples are classified into two classes:
\begin{enumerate}
\item Anosov or partially hyperbolic $\ZZ^n$-actions,
 and homogeneous actions of cocompact lattices related
 to Anosov or partially hyperbolic $\RR^n$-actions with $n \geq 2$
 ({\it e.g.} \cite{DK10,KL91,KS97,NT01}).
\item Isometric, or quasi-affine actions of lattices
 or groups with Property (T)
 ({\it e.g.} \cite{Be00,FM05,FM09,Zi85}).
\end{enumerate}
See Fisher's survey \cite{Fi07} for more related results.

One of the exceptions is an action of Baumslag-Solitar's group
 $BS(1,k)$ on the circle.
For $k \geq 2$, {\it Baumslag-Solitar's group $BS(1,k)$}
 is a finitely presented solvable group defined by
$BS(1,k)=\langle a,b \st aba^{-1}=b^k\rangle$.
It is isomorphic to
 a group generated by two affine transformations of the real line;
 $f(x)=kx$ and $g(x)=x+c$ with $c \neq 0$.
The natural extensions of $f$ and $g$ to $S^1=\RR \cup \{\infty\}$
 define a real analytic action $\rho_c$ of $BS(1,k)$ on $S^1$.
Remark that $\rho_c$ is conjugate to $\rho_1$ by
 a diffeomorphism $h(x)=c^{-1}x$.
\begin{thm}
[Burlsem and Wilkinson \cite{BW04}] 
\label{thm:BW}
Any real analytic action of $BS(1,k)$ on the circle is locally rigid.
In particular, the action $\rho_c$ is locally rigid.
\end{thm}
In the same paper,
 Burslem and Wilkinson also gave
 a smooth classification of $C^r$ actions of $BS(1,k)$ on $S^1$
 by using Navas' complete topological classification
 of $C^2$ solvable actions on one-dimensional manifolds (\cite{Na04}).
Guelman and Liousse \cite{GL11} extended
 the classification by Burslem and Wilkinson to $C^1$ actions
 by using Cantwell and Conlon's work \cite{CC02}
 on $C^1$ actions of $BS(1,k)$ on the circle
 or an closed interval, and Rivas' work \cite{Ri10}
 on $C^0$ action of $BS(1,k)$ on the real line.

Recently, some people have studied actions of
 Baumslag-Solitar like groups on higher dimensional manifolds.
McCarthy \cite{Mc10} proved the rigidity of trivial actions
 of a large class of abelian-by-cyclic groups on
 an arbitrary dimensional closed manifold.
Guelman and Liousse \cite{GL-pre} studied actions of $BS(1,k)$
 on surfaces, and gave a $C^\infty$ faithful action on the 2-torus
 which is not locally rigid even in topological sense.

In this paper, we study a natural higher dimensional analog
 of the standard $BS(1,k)$-action $\rho_c$.
For $n \geq 1$ and $k \geq 2$,
 we define a finitely generated solvable group $\Gamma_{k,n}$ by
\begin{equation*}
 \Gamma_{n,k}=\langle
 a,b_1,\dots,b_n \st ab_i a^{-1}=b_i^k,\; b_ib_j = b_j b_i
 \text{ for any }i,j=1,\dots,n
\rangle.
\end{equation*}
The group $\Gamma_{n,k}$ admits a natural action on
 the $n$-dimensional sphere $S^n$.
We identify $S^n$ with $\RR^n \cup \{\infty\}$
 by the stereographic projection.
For any basis $B=(v_1,\cdots,v_n)$ of $\RR^n$,
 define a $\Gamma_{n,k}$-action $\rho_B$ on $S^n$ by
\begin{itemize}
\item  $\rho_B^a(x) = kx$ and
 $\rho_B^{b_i}(x) = x+v_i$ for $x \in \RR^n=S^n \sm \{\infty\}$,
\item $\rho_B^a(\infty)=\rho_B^{b_i}(\infty)=\infty$.
\end{itemize}
The action $\rho_B$ preserves the standard conformal structure
 on $S^n$ and we call it {\it the standard action} associated to $B$.
For $n=1$ and $v_1=c \neq 0$, 
 the group $\Gamma_{1,k}$ is the Baumslag-Solitar group $BS(1,k)$
 and the action $\rho_B$ is the standard action $\rho_c$.
Therefore, $\rho_B$ is locally rigid by Theorem \ref{thm:BW} if $n=1$.
On the other hand,
 $\rho_B$ is not locally rigid for {\it any} basis $B$
 if $n \geq 2$ (see Proposition \ref{prop:classify}).
Hence, a direct analog of Theorem \ref{thm:BW} does not hold.

The aim of this paper is to show that the action $\rho_B$ exhibits
 a weak form of local rigidity for $n \geq 2$.

To state the main theorem,
 we recall basic concepts on rigidity of group actions.
Let $\Gamma$ be a discrete group
 and $G$ a topological group.
By $\Hom(\Gamma,G)$, we denote the set of
 homomorphism from $\Gamma$ to $G$.
For $\rho \in \Hom(\Gamma,G)$ and $\gamma \in \Gamma$,
 we put $\rho^\gamma=\rho(\gamma)$.
The set $\Hom(\Gamma,G)$ is naturally identified with
 a subset of a power set $G^\Gamma$.
The product topology on $G^\Gamma$ induces a topology on $\Hom(\Gamma,G)$.
When $G$ is Hausdorff,
 a sequence $(\rho_m)_{m \geq 1}$ in $\Hom(\Gamma,G)$ converges to $\rho$
 if and only if $\rho_m^\gamma$ converges to $\rho^\gamma$
 for any $\gamma \in \Gamma$.

Let $M$ be a smooth closed manifold.
In the below, all smooth maps and diffeomorphisms are of class $C^\infty$.
By $\Diff(M)$, we denote the group of diffeomorphisms of $M$.
It naturally becomes a topological group by the $C^\infty$-topology.
For a discrete group $\Gamma$,
 a smooth left $\Gamma$-action on $M$ is
 just a homomorphism from $\Gamma$ to $\Diff(M)$.
Hence, $\Hom(\Gamma,\Diff(M))$
 is identified with the space of (smooth left) $\Gamma$-actions on $M$.
We say that two actions $\rho_1 \in \Hom(\Gamma,\Diff(M_1))$
 and $\rho_2 \in \Hom(\Gamma,\Diff(M_2))$ are {\it smoothly conjugate}
 if there exists a diffeomorphism $h:M_1 \ra M_2$
 such that $\rho_2^\gamma \circ h = h \circ \rho_1^\gamma$
 for any $\gamma \in \Gamma$.
We also say that an action $\rho_0 \in \Hom(\Gamma,\Diff(M))$
 is {\it locally rigid} if
 there exists a neighborhood $\cU$ of $\rho_0$ in $\Hom(\Gamma,\Diff(M))$
 such that  any action $\rho$ in $\cU$ is smoothly conjugate to $\rho_0$.

Now, we are ready to state the main theorem of this paper.
\begin{Mthm}
Suppose $n,k \geq 2$.
Let $\rho_B$ be the standard $\Gamma_{n,k}$-action on $S^n$
 associated to a basis $B$ of $\RR^n$.
Then, there exists a neighborhood
 $\cU  \subset \Hom(\Gamma_{n,k},\Diff(S^n))$ of $\rho_B$
 such that any $\rho \in \cU$ is smoothly conjugate
 to $\rho_{B'}$ for some basis $B'=B'(\rho)$ of $\RR^n$.
In particular, any action in $\cU$ preserves a $C^\infty$
 conformal structure of $S^n$.
\end{Mthm}

The proof is divided into three steps:
 First, we show a local version of the main theorem,
 {\it i.e.}, rigidity of $\rho_B$ as a local action at $\infty$.
This is the main step of the proof.
Second, we prove that any perturbation of $\rho_B$ admits
 a global fixed point near $\infty$.
Finally, we extend the local conjugacy obtained in the first step
 to a global one.

The strategy for the first step
 is close to Burslem and Wilkinson's one in \cite{BW04}.
However, there is an essential difference from their case;
 the action $\rho_B$ admits non-trivial deformation.
The difficulty is that
 there seems no direct way
 to find a basis $B'=B'(\rho)$ such that $\rho$ is conjugate to $\rho_{B'}$
 for a given perturbation $\rho$ of $\rho_B$.
To overcome it, we follow Weil's idea in \cite{We64},
 where he controlled deformation of lattices of Lie groups
 by the first cohomology of a deformation complex.
Remark that Benveniste \cite{Be00} and Fisher \cite{Fi-pre}
 proved local rigidity of isometric actions
 by applying Weil's idea to $\Hom(\Gamma,\Diff(M))$.
In their cases, the deformation complex is infinite dimensional,
 and hence, they needed Hamilton's Implicit Function Theorem
 for tame maps between Fr\'echet spaces.
In our case, we reduce the deformation complex
 to a finite dimensional one
 and Weil's Implicit Function Theorem is sufficient.

In \cite{BW04}, Burslem and Wilkinson
 gave another proof of the first step above
 for $BS(1,k)$-actions on $S^1$.
They showed the existence of an invariant projective structure
 on a neighborhood of the global fixed point
 by using the Schwarzian derivative.
The author does not know
 whether there is an analogous proof for higher dimensional case.
Finding it seems an interesting problem.
 
\paragraph{Acknowledgements}
The author would like to thank an anonymous referee
 for valuable comments.

\section{Proof  of Main Theorem}

\subsection{Local version of the main theorem}
\label{sec:local}
Let $M_n(\RR)$ be the set of real square matrices
 of size $n$
 and $\GL_n(\RR)$ be the group of invertible matrices in $M_n(\RR)$.
We identify each element of $M_n(\RR)$ with
 an $n$-tuple of column vectors in $\RR^n$.
Under this identification,
 $\GL_n(\RR)$ is the set of bases of $\RR^n$.
By $\|\cdot\|$, we denote the Euclidean norm of $\RR^n$.
Let $\cS^r(\RR^n)$ be the set of symmetric $r$-multilinear maps
 from $(\RR^n)^r$ to $\RR^n$.
We define a norm $\|\cdot\|^{(r)}$ on $\cS^r(\RR^n)$ by
\begin{equation*}
 \|F\|^{(r)}
 = \sup\{\|F(\xi_1,\dots,\xi_r)\| \st \xi_1,\dots,\xi_r \in \RR^n,
 \|\xi_i\| \leq 1 \text{ for any }i\}.
\end{equation*}
Remark that
 $\|F(\xi_1,\dots,\xi_r)\| \leq \|F\| \cdot \|\xi_1\| \cdots 
 \|\xi_r\|$ for any $\xi_1, \dots,\xi_r \in \RR^n$
 and $\|A\|^{(1)}$ is the operator norm of $A \in M_n(\RR)=\cS^1(\RR^n)$.

Let $\cD(\RR^n,0)$ be the group of germs of
 local diffeomorphisms of $\RR^n$ at the origin.
For $F \in \cD(\RR^n,0)$,
 we denote the $r$-th derivative of $F$ at the origin by $D^{(r)}_0F$.
It is an element of $\cS^r(\RR^n)$.
For $r \geq 2$, we define {\it the $C^r_{loc}$-topology} on $\cD(\RR^n,0)$ by
 a pseudo-distance
 $d_{C^r_{loc}}(F,G)=\sum_{i=1}^r\|D^{(i)}_0F-D^{(i)}_0G\|^{(i)}$.
Remark that $d_{C^r_{loc}}$ is not a distance,
 and hence, the $C^r_{loc}$-topology is not Hausfdorff.

For a discrete group $\Gamma$,
 the $C^r_{loc}$-topology on $\Hom(\Gamma,\cD(\RR^n,0))$
 is naturally introduced as before.
We say that two local actions $P_1,P_2 \in \Hom(\Gamma,\cD(\RR^n,0))$
 are {\it smoothly conjugate}
 if there exists $H \in \cD(\RR^n,0)$ such that
 $P_2^\gamma \circ H = H \circ P_1^\gamma$ for any $\gamma \in \Gamma$.

Let $\bar{\phi}$ be a diffeomorphism from $S^n \sm \{0\}$
 to $\RR^n$ given by
\begin{equation*}
 \bar{\phi}(x)=\frac{1}{\|x\|^2}\cdot x.
\end{equation*}
For $B \in M_n(\RR)$, we define a local action
 $P_B \in \Hom(\Gamma_{n,k},\cD(\RR^n,0))$ by
\begin{equation*}
 P_B^\gamma=\bar{\phi} \circ \rho_B^\gamma \circ \bar{\phi}^{-1}. 
\end{equation*}
 
In this subsection, we prove 
 the following local version of the main theorem.
\begin{thm}
\label{thm:local}
 For $B \in \GL_n(\RR)$,
 there exists a $C^2_{loc}$-neighborhood $\cU$ of $P_B$
 such that any local action $P \in \cU$ is smoothly conjugate
 to $P_{B'}$ for some $B'=B'(P) \in \GL_n(\RR)$.
\end{thm}

The proof is divided into several steps.
First, we show the stability of linear part of $P^{b_i}$.
Let $\bar{F}$ be the element of $\cD(\RR^n,0)$ given by
\begin{equation*}
 \bar{F}(x)=k^{-1}x.
\end{equation*}
Notice that $P_B^a=\bar{F}$
 and $D^{(1)}_0 P_B^{b_i}=I$ for any $B \in M_n(\RR)$
 and $i=1,\dots,n$.
\begin{lemma}
\label{lemma:DG}
Let $m$ be a positive integer
 and $P_*$ be a local action in $\Hom(\Gamma_{n,k},\cD(\RR^m,0))$.
Suppose that $D^{(1)}_0P_*^a=k^{-1} I$ and $D^{(1)}_0P_*^{b_i}=I$
 for any $i=1,\dots,n$.
Then, there exists a $C^1_0$-neighborhood $\cU$ of $P_*$
 in $\Hom(\Gamma_{n,k},\cD(\RR^m,0))$
 such that $D_0^{(1)}P^{b_i}=I$
 for any $P \in \cU$ and $i=1,\dots,n$.
\end{lemma}
\begin{proof}
Put $c_{kj}=k!/[j!(k-j)!]$.
There exists $\delta>0$ such that
\begin{equation*}
\delta \cdot \left(k+k^2+\sum_{j=2}^kc_{kj}\delta^{j-2}(1+k\delta)\right)
 \leq \frac{1}{2}.
\end{equation*}
Take a $C^1_{loc}$-neighborhood $\cU$ of $P_*$ such that
 $\|D^{(1)}_0P^\gamma-D^{(1)}_0P_*^\gamma\| <\delta$
 for any $P \in \cU$ and $\gamma =a,b_1,\dots,b_m$.
Fix $P \in \cU$ and $i=1,\dots,n$.
We put $A=D^{(1)}_0P^a-D^{(1)}_0P_*^a$
 and $B=D^{(1)}_0P^{b_i}-D^{(1)}_0P_*^{b_i}$.
We need to show that $B=0$.
Since $D^{(1)}_0P_*^a=k^{-1}I$, $D^{(1)}_0P_*^{b_i}=I$,
 and $P^a \circ P^{b_i}=P^{b_i^k} \circ P^a$, we have
 $(k^{-1} I+A)(I+B)=(I+B)^k(k^{-1}I+A)$.
Hence,
\begin{align*}
(k-1)\|B\|^{(1)}
 & =\left\|
 kAB-k^2BA-\sum_{j=2}^kc_{kj}B^j(I+kA)
 \right\|^{(1)}\\
 & \leq
 \left(k\delta +k^2\delta+\sum_{j=2}^kc_{kj}\delta^{j-1}(1+k\delta))\right)
 \cdot \|B\|^{(1)}\\
 & \leq \frac{1}{2} \cdot \|B\|^{(1)}.
\end{align*}
Since $k \geq 2$, we obtain that $B=0$.
\end{proof}

Second, we show the stability of the linear part of $P^a$.
Let $\Pair{\cdot}{\cdot}$ be the Euclidean inner product of $\RR^n$.
For $v \in \RR^n$, we define $Q_v \in \cS^2(\RR^n)$ by
\begin{equation}
\label{eqn:D2Gv}
Q_v(\xi,\eta)
 =\Pair{\xi}{\eta} \cdot v-\Pair{\xi}{v} \cdot \eta
 -\Pair{\eta}{v} \cdot \xi.
\end{equation}
By a direct calculation, we can check that
\begin{equation*}
 D^{(2)}_0P_B^{b_i}=2Q_{v_i} 
\end{equation*}
 for any $B=(v_1,\dots,v_n) \in M_n(\RR)$ and $i=1,\dots,n$.
\begin{lemma}
\label{lemma:DF}
For any given $B \in \GL_n(\RR)$,
 there exists a $C^2_0$-neighborhood $\cU$ of $P_B$
 in $\Hom(\Gamma_{n,k},\cD(\RR^n,0))$
 such that $D_0^{(1)}P^a=k^{-1}I$
 for any $P \in \cU$.
\end{lemma}
\begin{proof}
For any $F,G \in \cD(\RR^n,0)$ with $D_0^{(1)}G=I$,
 it is easy to see that
\begin{align*}
D_0^{(2)}(F \circ G)
 & = D_0^{(2)}F +D_0^{(1)}F \circ D_0^{(2)}G,\\
D_0^{(2)}(G^k \circ F)
 & = D_0^{(2)}F+ k \cdot D_0^{(2)}G \circ (D_0^{(1)}F,D_0^{(1)}F).
\end{align*}

Put $B=(v_1,\dots,v_n)$.
Since $B=(v_1,\dots,v_n)$ is a basis of $\RR^n$,
 there exists a constant $\epsilon>0$ such that
 $\max_{i=1,\dots,n}\{\|A'v_i\|\} \geq \epsilon \|A'\|^{(1)} $
 for any $A' \in M_n(\RR)$.
By Lemma \ref{lemma:DG},
 there exists a $C^1_{loc}$-neighborhood $\cU_1$
 of $P_B$ such that $D_0^{(1)}P^{b_i}=I$  for any $P \in \cU_1$.
Let $\cU$ be a $C^2_{loc}$-open neighborhood of $P_B$
 consisting of $P \in \cU_1$ such that
\begin{equation*}
\max_{i=1,\dots,n}
\left\{3\|D_0^{(2)}P^{b_i}-2Q_{v_i}\|^{(2)}
 +\|k \cdot D_0^{(1)}P^a-I\|
 \cdot \|D^{(2)}_0 P^{b_i}\|^{(2)}\right\}
 < \epsilon.
\end{equation*}

Fix $P \in \cU_1$ and put
\begin{align*}
 A & = k \cdot D_0^{(1)}P^a-I,\\
 B_i & = D_0^{(2)}P^{b_i}-D_0^{(2)}P_B^{b_i}
 =D_0^{(2)}P^{b_i}-2Q_{v_i},\\
 C_i & = A \circ Q_{v_i} -2 Q_{v_i} \circ (A,I).
\end{align*}
We will show that $A=0$.
Since $P^a \circ P^{b_i}=P^{b_i^k} \circ P^a$, we have
\begin{equation*}
k^{-1}(I+A) \circ (2Q_{v_i}+B_i)
 = k \cdot (2Q_{v_i}+B_i) \circ (k^{-1}(I+A),k^{-1}(I+A)).
\end{equation*}
It implies that
\begin{align*}
2\|C_i\|^{(2)}
 & = \|A \circ B_i-2B_i \circ (A,I)
 -(2Q_{v_i}+B_i) \circ (A,A)\|^{(2)}\\
 & \leq \|A\|^{(1)} \cdot
 \left(3\|B_i\|^{(2)} +\|A\|^{(1)}\cdot \|D^{(2)}_0 P^{b_i}\|^{(2)}\right) \\
 & \leq  \epsilon \|A\|^{(1)}
\end{align*}
 for any $i=1,\dots,n$.
The definition of $Q_v$ also implies
 $C_i(v_i,v_i) =\|v_i\|^2 \cdot A v_i$,
 and hence, $\|C_i\|^{(2)} \geq \|A v_i\|$.
Therefore, we obtain
\begin{equation*}
 2\epsilon \|A\|^{(1)}
 \leq 2 \max_{i=1,\dots,n}\|A v_i\| \leq \epsilon \|A\|^{(1)}.
\end{equation*}
It implies that $A=0$, and hence, $D_0^{(1)}P^a= k^{-1} \cdot I$.
\end{proof}

Let $\cM_1'$ be the subset of $\Hom(\Gamma_{n,k},\cD(\RR^n,0))$
 consisting of local actions $P$ such that
 $P^a=\bar{F}$ and $D^{(1)}_0 P^{b_i}=I$  for any $i=1,\dots,n$.
Notice that $P_B$ is an element of $\cM_1'$ for any $B \in \GL_n(\RR)$.
\begin{prop}
\label{prop:reduced} 
Let $B$ be an element of $\GL_n(\RR)$.
For any given $C^2_{loc}$-neighborhood $\cU_0$ of $P_B$
 in $\Hom(\Gamma_{n,k},\cD(\RR^n,0))$,
 there exists another $C^2_{loc}$-neighborhood $\cU$ of $P_B$
 such that
 any $P \in \cU$ is smoothly conjugate to
 a local action in $\cU_0 \cap \cM_1'$.
\end{prop}
\begin{proof}
By Lemmas \ref{lemma:DG} and \ref{lemma:DF},
 there exists a $C^2_{loc}$-neighborhood $\cU_1$ of $P_B$
 in $\Hom(\Gamma_{n,k},\cD(\RR^n,0))$ such that
 $D^{(1)}_0 P^a= k^{-1} \cdot I$ and $D^{(1)}_0P^{b_i}=I$ 
 for any $P \in \cU_1$ and $i=1,\dots,n$.
Fix $P \in \cU_1$.
It is known that if a local diffeomorphism $F \in \cD(\RR^n,0)$
 satisfies  $D^{(1)}_0F=\alpha I$ for some $0<\alpha<1$
 then it is smoothly linearizable
 (see {\it e.g.} \cite[Theorem 6.6.6]{KH}).
Hence, their exists $H \in \cD(\RR^n,0)$ such that
 $D^{(1)}_0 H=I$ and $\bar{F} =H \circ P^a \circ H^{-1}$.
We define a local action $P^H \in \Hom(\Gamma_{n,k},\cD(\RR^n,0))$ by
 $(P^H)^\gamma=H \circ P^\gamma \circ H^{-1}$.
Since $D^{(1)}_0 (P^H)^{b_i}=D^{(1)}_0P^{b_i}=I$,
 the local action $P^H$ is contained in $\cM_1'$.
From the equation $D^{(2)}_0(H \circ \bar{F})=D^{(2)}_0(P^a \circ H)$,
 we obtain $(k-1)D^{(2)}_0H=k^2 D^{(2)}_0 P^a$.
Hence, there exists
 a small $C^2_{loc}$-neighborhood $\cU \subset \cU_1$ of $P_B$
 such that  $P^H \in \cU_0$ for any $P \in \cU$.
\end{proof}

Following Weil's idea,
 we reduce Theorem \ref{thm:local} to exactness of a linear complex.
Put
\begin{align*}
\cM_0 & = \GL_n(\RR) \times \GL_n(\RR),\\
\cM_1 & = \left\{(G_i)_{1 \leq i \leq n} \in \cD(\RR^n,0)^n
 \;\left|\;
 D^{(1)}_0G_i=I, \bar{F} \circ G_i = G_i^k \circ \bar{F}
 \right.
 \text { for any } i \right\},\\
\cM_2 & =\{(C_{ij})_{1 \leq i<j \leq n} \st C_{ij} \in \cS^3(\RR^n)\}
 = (\cS^3(\RR^n))^{n(n-1)/2}.
\end{align*}
Define maps $\Phi:\cM_0 \ra \cM_1$ and $\Psi:\cM_1 \ra \cM_2$ by
\begin{align*}
\Phi(A,B) & = (A \circ P_B^{b_i} \circ A^{-1})_{1 \leq i \leq n},\\
\Psi((G_i)_{1 \leq i \leq n})
 & = \left(\frac{1}{4}
  \left[D^{(3)}_0(G_i \circ G_j)-D^{(3)}_0(G_j \circ G_i)\right]
 \right)_{1 \leq i<j \leq n}
\end{align*}
By $O_{\cM_2}$, we denote the zero element of $\cM_2=S^3(\RR^n)^{n(n-1)/2}$.
Then,
\begin{equation*}
 \Psi \circ \Phi(A,B) = O_{\cM_2},\hsp
 \Psi(P^{b_1},\dots,P^{b_n})=O_{\cM_2} 
\end{equation*}
 for any $(A,B) \in \cM_0$ and $P \in \cM_1'$.
Moreover, if $\Phi(A,B)=(P^{b_1},\dots,P^{b_n})$,
 then $P$ is smoothly conjugate to $P_B$ by the linear map $A$.

The following is a direct corollary of Proposition \ref{prop:reduced}.
\begin{cor}
\label{cor:non-linear exact}
To prove Theorem \ref{thm:local},
 it is sufficient to show the existence of
 a $C^2_{loc}$-neighborhood $\cV_*$ of $(P_B^{b_i})_{1 \leq i \leq n}$
 in $\cM_1$ such that
\begin{equation*}
\Psi^{-1}(O_{\cM_2}) \cap \cV_*
 =\Img \Phi \cap \cV_* 
\end{equation*}
 for any given $B \in \GL_n(\RR)$.
\qed
\end{cor}

Let us recall Weil's Implicit Function Theorem.
\begin{thm}
[Weil, \cite{We64}] 
\label{thm:Weil}
Let $\Phi_0:M_0 \ra M_1$ and $\Phi_1:M_1 \ra M_2$
 be smooth maps between manifolds $M_0$, $M_1$, and $M_2$.
Suppose that $\Phi_1 \circ \Phi_0$ is a constant map
 with value $x_2 \in M_2$.
If $\Ker (D\Phi_1)_{x_1}=\Img (D\Phi_0)_{x_0}$
 for $x_0 \in M_0$ and $x_1=\Phi_0(x_0) \in M_1$,
 then there exists a neighborhood $U$ of $x_1$ such that 
 $\Img \Phi_0 \cap U = \Phi_1^{-1}(x_2) \cap U$.
\end{thm}

The spaces $\cM_0$ admits a natural smooth structure
 as an open subset of a finite dimensional vector space $M_n(\RR)^2$.
The space $\cM_2=(\cS^3(\RR^n))^{n(n-1)/2}$ also does
 as a finite dimensional vector space.
If the maps $\Phi$ and $\Psi$ are smooth
 with respect to some smooth structure on $\cM_1$ compatible to
 the $C^2_{loc}$-topology
 and they satisfy
 $\Ker D\Psi_{(P_B^{b_1},\dots,P_B^{b_n})}=\Img D\Phi_{(I,B)}$,
 then Theorem \ref{thm:local} follows from
 Corollary \ref{cor:non-linear exact} and Weil's theorem.
To introduce a smooth structure on $\cM_1$,
 we define a map $\Theta:\cM_1 \ra \cS^2(\RR^n)^n$ by
\begin{equation*}
\Theta(G_1,\dots,G_n) = \frac{1}{2}(D^{(2)}_0G_1,\dots,D^{(2)}_0 G_n).
\end{equation*}
\begin{lemma}
\label{lemma:Theta} 
 The map $\Theta$ is a homeomorphism
 with respect to the $C^2_{loc}$-topology on $\cM_1$.
\end{lemma}
\begin{proof}
Since $D^{(1)}_0 G_i=I$ for any $(G_1,\dots,G_n) \in \cM_1$
 and any $i$, the map $\Theta$ is continuous
 by the definition of the $C^2_{loc}$-topology.

Next, we show that $\Theta$ is surjective.
Put $(e_1,\dots,e_n)=I$ and take $Q \in \cS^2(\RR^n)$.
Let $G_Q^t \in \cD(\RR^n,0)$ be the time-$t$ map of the local flow
 generated by the quadratic vector field $X_Q(x)=Q(x,x)$.
Then, $G_Q^t(0)=0$, $D^{(1)}_0G_Q^t=I$,
 and $\bar{F} \circ G_Q^t =G_Q^{kt} \circ \bar{F}$ for any $t \in \RR$.
Since
\begin{align*}
 \left.\frac{d}{dt}[D^{(2)}_0G_Q^t](e_i,e_j)\right|_{t =t_0}
& =\left. \frac{\del}{\del t}\frac{\del^2}{\del x_i \del x_j} G_Q^t(x)
  \right|_{(x,t)=(0,t_0)}\\
& =\left. \frac{\del^2}{\del x_i \del x_j}\frac{\del}{\del t} G_Q^t(x)
  \right|_{(x,t)=(0,t_0)}\\
& =\left. \frac{\del^2}{\del x_i \del x_j}
  Q(G_Q^{t_0}(x),G_Q^{t_0}(x)) \right|_{x=0}\\
& =2Q(D^{(1)}_0G_Q^{t_0}(e_i),D^{(1)}_0G_Q^{t_0}(e_j))\\
& = 2Q(e_i,e_j)
\end{align*}
 for any $i,j=1,\dots,n$ and $t_0 \in \RR$,
 we have $D^{(2)}_0G_Q^t=2tQ$ for any $t$.
Therefore, $\Theta(G_{Q_1}^1,\dots,G_{Q_n}^1)=(Q_1,\dots,Q_n)$
 for any $(Q_1,\dots,Q_n) \in \cS^2(\RR^n)^n$.

Finally, we show that $\Theta$ is injective.
Remark that the bijectivity of $\Theta$ implies that it is an open map.
Take $G_1,G_2 \in \cD(\RR^N,0)$ such that
 that $D_0^{(1)}G_i=I$ and $\bar{F} \circ G_i = G_i^k \circ \bar{F}$
 for $i=1,2$, and $D_0^{(2)}G_1=D_0^{(2)}G_2$.
We will show that $G_1=G_2$.
For $R>0$, we put $B_R=\{z \in \RR^N \st \|z\| \leq R\}$.
Fix representatives $\tilde{G}_i$ of $G_i$ for each $i=1,2$.
Since $G_i(0)=0$ and $D_0^{(1)}G_i=I$,
 there exists $R_0>0$ and $1<c<\sqrt[4]{k}$ such that
\begin{itemize}
\item $\tilde{G}_2^m \circ \tilde{G}_1^{m'}$ is well-defined on $B_{R_0}$
 for any $m, m'=1,\dots,k$,
\item $\bar{F} \circ \tilde{G}_i=\tilde{G}_i^k \circ \bar{F}$
 on $B_{R_0}$ for $i=1,2$,
\item $\max\{\|\tilde{G}_1^m(z)\|
 ,\|\tilde{G}_2 \circ \tilde{G}_1^m(z)\|\} \leq c\|z\|$,
 and $\|\tilde{G}_2^m(z)-\tilde{G}_2^m(z')\| \leq c\|z-z'\|$
 for any $z,z' \in B_{R_0}$ and $m=1,\dots,k$.
\end{itemize}
For $0<R \leq R_0$, we put
\begin{equation*}
 \Delta(R)=\sup_{z \in B_R}\frac{\|\tilde{G}_1(z)-\tilde{G}_2(z)\|}{\|z\|^3}.
\end{equation*}
Since $D_0^{(2)}G_1=D_0^{(2)}G_2$, then $\tilde{G}_1-\tilde{G}_2$
 is of at least third order at the origin.
Hence, $\Delta(R)$ is finite.
For any $z \in B_{R_0}$ and $m=1,\dots,k$,
 we have $\max\{
 \|\tilde{G}_1^m(k^{-1}z)\|, \|\tilde{G}_2 \circ \tilde{G}_1^m(k^{-1}z)\|
 \} \leq (c/k)\|z\| \leq R_0$,
 and hence,
\begin{align*}
\|\tilde{G}_1(z)-\tilde{G}_2(z)\|
 & = k \cdot \|\bar{F} \circ \tilde{G}_1(z)-\bar{F} \circ \tilde{G}_2(z)\|\\
 & =k  \cdot \|\tilde{G}_1^k(k^{-1}z)-\tilde{G}_2^k(k^{-1}z)\|\\
 & \leq k \sum_{m=1}^k
 \|\tilde{G}_2^{m-1} \circ \tilde{G}_1^{k-m+1}(k^{-1}z)
 -\tilde{G}_2^m \circ \tilde{G}_1^{k-m}(k^{-1}z)\|\\
 & \leq kc\sum_{m=1}^k
 \|\tilde{G}_1^{k-m+1}(k^{-1}z)
 -\tilde{G}_2 \circ \tilde{G}_1^{k-m}(k^{-1}z)\|.
\end{align*}
Since $\|\tilde{G}_1^{k-m}(k^{-1}z)\| \leq (c/k)\|z\| \leq  R_0$,
 it implies that
\begin{align*}
\|\tilde{G}_1(z)-\tilde{G}_2(z)\|
 & \leq k^2c \cdot \Delta(R_0) \cdot [(c/k) \cdot\|z\|]^3
 = (c^4/k) \cdot \Delta(R_0) \cdot \|z\|^3.
\end{align*}
Therefore, $\Delta(R_0) \leq (c^4/k)\Delta(R_0)$.
Since $c<\sqrt[4]{k}$, we have $\Delta(R_0)=0$,
 and hence, $\tilde{G}_1=\tilde{G}_2$ on $B_{R_0}$
\end{proof}

Since
\begin{equation*}
  D^{(2)}_0(A \circ P_B^{b_i} \circ A^{-1})=
  A \circ D^{(2)}_0 P_B^{b_i} \circ (A^{-1},A^{-1})
 = 2A \circ Q_{v_i} \circ (A^{-1},A^{-1}).
\end{equation*}
 for $A \in \GL_n(\RR)$ and $B=(v_1,\dots,v_n) \in \GL_n(\RR)$,
 the map $\Theta \circ \Phi$ satisfies
\begin{equation}
\label{eqn:Phi}
 (\Theta \circ \Phi)(A,B)
 =(A \circ Q_{v_i} \circ (A^{-1},A^{-1}))_{1 \leq i \leq n}.
\end{equation}
Hence, $\Theta \circ \Phi$ is smooth.
For $Q,Q' \in \cS^2(\RR^n)$,
 we define the bracket $[Q,Q'] \in \cS^3(\RR^n)$ by
\begin{align*}
 [Q,Q'](\xi,\eta,\theta)
 & =\left\{Q(\xi,Q'(\eta,\theta))+Q(\eta,Q'(\theta,\xi))
 +Q(\theta,Q'(\xi,\eta))\right\}\\
 & \quad -\left\{Q'(\xi,Q(\eta,\theta))+Q'(\eta,Q(\theta,\xi))
 +Q'(\theta,Q(\xi,\eta))\right\}
\end{align*}
It can be checked that
\begin{equation}
 [D_0^{(2)}G_1,D_0^{(2)}G_2]
 = D_0^{(3)}(G_1 \circ G_2) -D_0^{(3)}(G_2 \circ G_1)
\end{equation}
 for any $G_1,G_2 \in \cD(\RR^N,0)$ with $D_0^{(1)}G_1=D_0^{(1)}G_2=I$.
Therefore,
\begin{equation}
\label{eqn:Psi}
 (\Psi \circ \Theta^{-1})(Q_1,\dots,Q_n)=([Q_i,Q_j])_{1 \leq i<j \leq n}.
\end{equation}
Since the bracket is bi-linear,
 the map $\Psi \circ \Theta^{-1}$ is a smooth map.

For $B=(v_1,\dots,v_n) \in \GL_n(\RR)$, we put
\begin{align*}
L^\Phi_B & =D(\Theta \circ \Phi)_{(I,B)},\\
L^\Psi_B & =D(\Psi \circ \Theta^{-1})_{(Q_{v_1},\dots,Q_{v_n})}.
\end{align*}
We identify the tangent spaces of $\cM_0$ and $\cS^2(\RR^n)^n$
 of each point with $M_n(\RR^n)^2$ and $\cS^2(\RR^n)^n$, respectively.
Then, Equations (\ref{eqn:Phi}) and (\ref{eqn:Psi}) imply that
\begin{align*}
L^\Phi_B(A',B') &
 = (A' \circ Q_{v_i} - Q_{v_i} \circ (A',I)-Q_{v_i} \circ (I,A')
  + Q_{\omega_i})_{1 \leq i \leq n}\\
L^\Psi_B(q_1,\dots,q_n) &
 = ([q_i,Q_{v_j}]-[q_j,Q_{v_i}])_{1 \leq i <j \leq n}
\end{align*}
 for any $(A',B') \in M_n(\RR)^2$
 with $B'=(\omega_1,\dots,\omega_n)$
 and any $(q_1,\dots,q_n) \in \cS^2(\RR^n)^n$.
The following proposition can be shown by a formal computation
 and we postpone the proof until Section \ref{sec:vanishing}.
\begin{prop}
\label{prop:vanishing} 
$\Ker L^\Psi_B=\Img L^\Phi_B$.
\end{prop}
Theorem \ref{thm:local} follows from
 Corollary \ref{cor:non-linear exact},
 Theorem \ref{thm:Weil}, and the proposition
 since $H$ is a homeomorphism between $\cM_1$ and $\cS^2(\RR^n)^n$.

\subsection{From local to global}
In this subsection, we prove the main theorem.
For a discrete group $\Gamma$
 and a $\Gamma$-action $\rho$ on a manifold $M$,
 we say that a point $p \in M$ is {\it a global fixed point}
 if $\rho^\gamma(p)=p$ for any $\gamma \in \Gamma$.
Remark that the point $\infty$ is the unique global fixed point of $\rho_B$
 for any $B \in \GL_n(\RR)$.

{\it In this subsection, we assume that $n \geq 2$}
 since the case $n=1$ was already shown by Burslem and Wilkinson.
First, we show that any local conjugacy
 to the standard $\Gamma_{n,k}$-action extends to a global one.
\begin{prop}
\label{prop:global conj}
Suppose that an action $\rho \in \Hom(\Gamma_{n,k},\Diff(M))$
 admits a global fixed point $p_\infty$
 and there exists a smooth coordinate $\phi$ of $S^n$ at $p_\infty$
 and $B \in \GL_n(\RR)$
 such that $\phi(p_\infty)=0$
 and $\phi \circ \rho^\gamma \circ \phi^{-1}=P_B^\gamma$
 as elements of $\cD(\RR^n,0)$ for any $\gamma \in \Gamma_{n,k}$.
Then, $\rho$ is smoothly conjugate to $\rho_B$.
\end{prop}
\begin{proof}
Recall that $\bar{\phi}:S^n \ra \RR^n$ is the local coordinate at $\infty$
 given by $\bar{\phi}(x)=(1/\|x\|^2) \cdot x$
 and the local action $P_B$
 is defined by $P_B^\gamma = \bar{\phi} \circ \rho_B^\gamma \circ \bar{\phi}$.
We put $U_r=S^n \sm [-r,r]^n$ for $r>0$
 and $\Lambda_b=\{b_1^{\pm 1}, \dots, b_n^{\pm 1}\}$.
By assumption,
 there exists $R>0$ and a neighborhood $U'$ of $p_\infty$ such that
\begin{equation*}
\phi \circ \rho^{\gamma} \circ \phi^{-1}
 =\bar{\phi} \circ \rho_B^\gamma \circ \bar{\phi}^{-1}
\end{equation*}
 on $\bar{\phi}(U_R)$
 for any $\gamma \in \{a^{\pm 1}\} \cup \Lambda_b$.
Since $\rho_B^{b_1^m}(x)$ converges to $\infty$
 as $n$ goes to infinity for any $x \in S^n$,
 we can take $m_x \geq 0$ such that $\rho_B^{b_1^{m_x}}(x)$ is
 contained in $U_R$.
Define a map $h:S^n \ra S^n$ by
\begin{equation*}
 h(x)=\rho^{b_1^{-m_x}} \circ (\phi^{-1} \circ \bar{\phi})
 \circ \rho_B^{b_1^{m_x}}(x)
\end{equation*}
First, we see that $h(x)$ does not depend on the choice of $m_x$.
Suppose that $\rho_B^{b_1^m}(x)$ is contained in $U_R$.
Since $\rho_B^\gamma$ is a translation for any $\gamma \in \Lambda_b$
 and $S^n \sm U_R=[-R,R]^n$ is a convex subset of $\RR^n$,
 there exists a sequence $(\gamma_j)_{1 \leq j \leq l}$ 
 of elements of $\Lambda_b$ such that
 $b_1^m=\gamma_l \cdots \gamma_1 b_1^{m_x}$
 and $\rho_B^{\gamma_j\cdots\gamma_1 b_1^{m_x}}(x)$ is contained in $U_R$
 for any $j=1,\dots,l$.\footnote{We need $n \geq 2$ here.}
Then, 
\begin{equation*}
 \rho^{\gamma_{j+1}} \circ (\phi^{-1} \circ \bar{\phi})
 \circ \rho_B^{\gamma_j\cdots\gamma_1 b_1^{m_x}}(x)
 = (\phi^{-1} \circ \bar{\phi}) \circ
 \rho_B^{\gamma_{j+1}\gamma_j\cdots\gamma_1 b_1^{m_x}}(x).
\end{equation*}
This implies that
\begin{align*}
 \rho^{b_1^{-m}} \circ (\phi^{-1} \circ \bar{\phi}) \circ \rho^{b_1^m}(x)
 & =\rho^{b_1^{-m}}
 \circ (\phi^{-1} \circ \bar{\phi}) \circ \rho_B^{\gamma_l \cdots \gamma_1} 
 \circ \rho^{b_1^{m_x}}(x)\\
 & =\rho^{b_1^{-m}} \circ \rho^{\gamma_l \cdots \gamma_1} 
 \circ (\phi^{-1} \circ \bar{\phi}) \circ \rho^{b_1^{m_x}}(x)\\
 & =\rho^{b_1^{-m_x}}
 \circ (\phi^{-1} \circ \bar{\phi}) \circ \rho^{b_1^{m_x}}(x)\\
 & =h(x).
\end{align*}
Therefore, $h(x)$ does not depend on the choice of $m_x$.

For any given $x_0 \in S^n$,
 there is a choice of $(m_x)_{x \in S^n}$ which is constant on
 a small neighborhood of $x_0$.
This implies that $h$ is a locally diffeomorphic at $x_0$,
 and hence, $h$ is a covering map.
Since $S^n$ is simply-connected, $h$ is diffeomorphism.

It is easy to see that $h \circ \rho_B^\gamma=\rho^\gamma \circ h$
 for any $\gamma \in \Lambda_b$.
For any given $x \in S^n$,
 there exists $m \geq 1$ such that $\rho_B^{b_1^{km}}(x)$
 is contained in $U_R$.
Then,
\begin{align*}
h \circ \rho_B^a(x)
 & = \rho^{b_1^{-km}} \circ (\phi^{-1} \circ \bar{\phi})
  \circ \rho_B^{b_1^{km}} \circ \rho_B^a(x)\\
 & = \rho^{b_1^{-km}} \circ (\phi^{-1} \circ \bar{\phi})
   \circ \rho_B^a \circ \rho_B^{b_1^{m}}(x)\\
 & = \rho^{b_1^{-km}}\circ \rho^a
   \circ (\phi^{-1} \circ \bar{\phi}) \circ \rho_B^{b_1^{m}}(x)\\
 & = \rho^a \circ \rho^{b_1^{-m}}
   \circ (\phi^{-1} \circ \bar{\phi}) \circ \rho_B^{b_1^{m}}(x)\\
 & = \rho^a \circ h(x).
\end{align*}
Therefore, $h$ is a smooth conjugacy between $\rho_B$ and $\rho$.
\end{proof}

Next, we give a criterion
 for the persistence of a global fixed point of a $\Gamma_{n,k}$-action.
\begin{lemma}
\label{lemma:fixed point}
Let $M$ be a manifold and
 $\rho$ be an action in $\Hom(\Gamma_{n,k},\Diff(M))$.
Suppose that $\rho_0$ has a global fixed point $p_0$ such that
 $(D\rho_0^a)_{p_0}=k^{-1}I$ and $(D\rho_0^{b_i})_{p_0}=I$ for any
 $i=1,\dots,n$.
Then, there exists a neighborhood
 $\cU \subset \Hom(\Gamma_{n,k},\Diff(M))$ of $\rho_0$
 and a continuous map $\hat{p}:\cU \ra M$ 
 such that $\hat{p}(\rho_0)=p_0$
 and  $\hat{p}(\rho)$ is a global fixed point of $\rho$
 for any $\rho \in \cU$.
\end{lemma}
\begin{proof}
Take $k^{-1}<\lambda<1$ and
 $\delta>0$ so that $\lambda+k\delta<1$.
Fix an open neighborhood $U$ of $p_0$ and 
 a local coordinate $\phi:U \ra \RR^n$.
There exist convex neighborhoods $V$ and $V_1$ of $\phi(p_0)$
 and a neighborhood $\cU_0$ of $\rho_0$
 which satisfy the following conditions
 for any $\rho \in \cU_0$ and $i=1,\dots,n$;
\begin{itemize}
\item $\phi \circ \rho^{a^l b_i^m} \circ \phi^{-1}$
 is well-defined on $V$ for any $l=0,1$ and $m=0,\dots,k$.
\item $\phi \circ \rho^{b_i} \circ \phi^{-1}(V_1) \subset V$.
\item $\|D(\phi \circ \rho^a \circ \phi^{-1})_z\|<\lambda$
 and $\|D(\phi \circ \rho^{b_i^m} \circ \phi^{-1})_z-I\|<\delta$
 for any $z \in V$ and $m=1,\dots,k$.
\end{itemize}
By the persistence of attracting fixed point,
 there exists a neighborhood
 $\cU \subset \cU_0$ of $\rho_0$
 and a continuous map $\hat{p}:\cU \ra \phi^{-1}(V_1 \cap V)$ 
 such that $\hat{p}(\rho_0)=p_0$ and
  $\hat{p}(\rho)$ is an attracting fixed point of $\rho^a$
 for any $\rho \in \cU$.
Since $\rho_0^{b_i}(p_0)=p_0$,
 by replacing $\cU$ with a smaller neighborhood of $\rho_0$,
 we may assume that $\rho^{b_i}(\hat{p}(\rho)) \in \phi^{-1}(V_1 \cap V)$
 for any $\rho \in \cU$ and $i=1,\dots,n$.

Fix $i=1,\dots,n$ and $\rho \in \cU$.
Put $z_*=\phi(\hat{p}(\rho))$,
 $F=\phi \circ \rho^a \circ \phi^{-1}$,
 and $G=\phi \circ \rho^{b_i} \circ \phi^{-1}$.
We will show $G(z_*)=z_*$.
Since $z_*$ and $G(z_*)$ are contained in $V$,
\begin{gather*}
\|F \circ G(z_*)-F(z_*)\| \leq \lambda \|G(z_*)-z_*\|,\\
 \|(G^{m+1}(z_*)-G(z_*))-(G^m(z_*)-z_*)\| \leq \delta \|G(z_*)-z_*\|.
\end{gather*}
 for $m=0,\dots,k-1$.
Since $F \circ G=G^k \circ F$ and $F(z_*)=z_*$,
 the former implies
\begin{equation*}
 \|G^k(z_*)-z_*\| \leq \lambda \|G(z_*)-z_*\|.
\end{equation*}
Hence, 
\begin{align*}
k \cdot \|z_*-G(z_*)\|
 & \leq \|G^k(z_*)-z_*\|+\sum_{m=0}^{k-1}
 \|G^{m+1}(z_*)-G^m(z_*)-G(z_*)+z_*\| \\
 & \leq (\lambda+k\delta) \cdot \|G(z_*)-z_*\|.
\end{align*}
Since $\lambda+k\delta<1$, this implies $G(z_*)=z_*$.
Therefore, $\hat{p}(\rho)$ is a global fixed point of $\rho$.
\end{proof}

Now, we prove the main theorem.
\begin{proof}
[Proof of Main Theorem] 
Take open neighborhoods $U \subset S^n$ of $\infty$
 and $V \subset \RR^n$ of $0$,
 and a family $(\phi_p)_{p \in U}$ of diffeomorphisms from $U$ to $V$
 such that $\phi_\infty=\bar{\phi}$,
 $\phi_p(p)=0$ for any $p \in U$,
 and the map $(p,q) \mapsto \phi_p(q)$ is smooth.
Fix $B \in \GL_n(\RR)$.
The action $\rho_B$ satisfies the assumption
 of Lemma \ref{lemma:fixed point}.
Hence, there exists a neighborhood $\cU_1$ of $\rho_B$
 and a continuous map $\hat{p}:\cU_1 \ra U$
 such that $\hat{p}(\rho)$ is a global fixed point of $\rho$
 for any $\rho \in \cU_1$.
We define a local action $P_\rho \in \Hom(\Gamma_{n,k},\cD(\RR^n,0))$
 by $P_\rho^\gamma
 =\phi_{\hat{p}(\rho)} \circ \rho^\gamma \circ \phi_{\hat{p}(\rho)}^{-1}$.
Then, the map $\rho \mapsto P_\rho$ is $C^2_{loc}$-continuous map
 from $\cU_1$ to $\Hom(\Gamma_{n,k},\cD(\RR^n,0))$.
By Theorem \ref{thm:local},
 there exists a neighborhood $\cU \subset \cU_1$ of $\rho_B$
 such that $P_\rho$ is smoothly conjugate to $P_{B'}$
 for some $B'=B'(\rho) \in \GL_n(\RR)$
 for any $\rho \in \cU$.
By Proposition \ref{prop:global conj},
 $\rho$ is smoothly conjugate to $\rho_{B'}$.
\end{proof}

\subsection{Proof of Proposition \ref{prop:vanishing}}
\label{sec:vanishing}
In this subsection,
 we give a proof of the following proposition,
 which we have not shown in Subsection \ref{sec:local}.
\begin{propV}
$\Ker L^\Psi_B=\Img L^\Phi_B$. 
\end{propV}
Our proof is formal and lengthy computation.
It may be interesting to find a more geometric proof.

Fix $B = (v_1,\dots,v_n) \in \GL_n(\RR)$.
Recall that the linear maps
 $L^\Phi_B:M_n(\RR^n)^2 \ra \cS^2(\RR^n)^n$
 and $L^\Psi_B:\cS^2(\RR^n)^n \ra \cS^3(\RR^n)^{n(n-1)/2}$
 are given by
\begin{align*}
L^\Phi_B(A',B') &
 = (A' \circ Q_{v_i} - Q_{v_i} \circ (A',I)-Q_{v_i} \circ (I,A')
 + Q_{\omega_i})_{1 \leq i \leq n}\\
L^\Psi_B(q_1,\dots,q_n) &
 = ([q_i,Q_{v_j}]-[q_j,Q_{v_i}])_{1 \leq i <j \leq n}
\end{align*}
 for any $(A',B') \in M_n(\RR)^2$
 with $B'=(\omega_1,\dots,\omega_n)$
 and any $(q_1,\dots,q_n) \in \cS^2(\RR^n)^n$,
 where
\begin{equation}
\label{eqn:D2Gv 2}
Q_v(\xi,\eta)
 =\Pair{\xi}{\eta} \cdot v-\Pair{\xi}{v} \cdot \eta
 -\Pair{\eta}{v} \cdot \xi.
\end{equation}
 and
\begin{align*}
 [Q,Q'](\xi,\eta,\theta)
 & =\left\{Q(\xi,Q'(\eta,\theta))+Q(\eta,Q'(\theta,\xi))
 +Q(\theta,Q'(\xi,\eta))\right\}\\
 & \quad -\left\{Q'(\xi,Q(\eta,\theta))+Q'(\eta,Q(\theta,\xi))
 +Q'(\theta,Q(\xi,\eta))\right\}.
\end{align*}

First, we reduce the problem to the case $B=I$.
\begin{lemma}
\label{lemma:standard}
For $B,B' \in \GL_n(\RR)$,
 $\Ker L_B^\Psi=\Img L_B^\Phi$
 if and only if $\Ker L_{B'}^\Psi=\Img L_{B'}^\Phi$.
\end{lemma}
\begin{proof}
Put $B=(v_1,\dots,v_n)$ and $B'=(w_1,\dots,w_n)$.
Take $A=(a_{ij}) \in \GL_n(\RR)$ such that $B'=BA$.
Since the map $v \mapsto Q_v$ is linear,
\begin{equation*}
 (Q_{w_1},\dots,Q_{w_n})=(Q_{v_1},\dots,Q_{v_n}) \cdot A. 
\end{equation*}
It implies that
$\Img L_{B'}^\Phi=\Img L_{B}^\Phi \cdot A$.
For $(q_1,\dots,q_n) \in \Ker L_B^\Psi$ we have
\begin{align*}
\left[\left(\sum_{k=1}^n a_{ki}q_k\right), Q_{w_j}\right]
 -\left[\left(\sum_{l=1}^n a_{lj}q_l\right), Q_{w_i}\right]
 & =\sum_{k,l=1}^n a_{ki}a_{lj}([q_k,Q_{v_l}]-[q_l,Q_{v_k}])
 =0.
\end{align*}
Hence, $\Ker L_{B}^\Psi \cdot A$ is a subspace of $\Ker L_{B'}^\Psi$.
Similarly, $\Ker L_{B'}^\Psi \cdot A^{-1}$ is a subspace of $\Ker L_{B}^\Psi$.
Therefore, $\Ker L_{B'}^\Psi=\Ker L_{B}^\Psi \cdot A$.
\end{proof}

By the lemma, it is sufficient
 to show Proposition \ref{prop:vanishing} for $B=I$.
Put $I=(e_1,\dots,e_n)$.
It is easy to check the following properties of $Q_v$.
\begin{lemma}
\label{lemma:inj 0} 
For $v \in \RR^n$ and mutually disjoint $i,j,k=1,\dots,n$,
\begin{gather*}
 Q_{e_i}(e_i,v)=Q_{e_i}(v,e_i)=-v,\\
 Q_{e_i}(e_j,e_j)=e_i,\\
 Q_{e_i}(e_j,e_k)=0.
\end{gather*}
\qed
\end{lemma}

Let $W$ be the subspace of $\cS^2(\RR^n)^n$ consisting of
 $(q_1,\dots,q_n)$ such that
\begin{gather}
\label{eqn:H 1}
q_j(e_j,e_j)=0, \\
\label{eqn:H 2}
\Pair{e_i}{q_j(e_i,e_i)}+\Pair{e_j}{q_i(e_j,e_j)}=0,\\
\label{eqn:H 3}
\Pair{e_1}{q_1(e_j,e_j)}=0
\end{gather}
 for any $i,j=1,\dots,n$.

\begin{lemma}
\label{lemma:H}
If $\Ker L_I^\Psi \cap W= \{0\}$,
 then $\Ker L_I^\Psi =\Img L_I^\Phi$.
\end{lemma}
\begin{proof}
We show that $\cS^2(\RR^n)^n=W+\Img L_I^\Phi$.
Once it is shown,
 then the assumption $\Ker L_I^\Psi \cap W= \{0\}$
 implies $\Ker L_I^\Psi =\Img L_I^\Phi$
 since $\Img L_I^\Phi \subset \Ker L_I^\Psi$. 

For $A,B \in M_n(\RR)$,
 let  $q^{A,B}_j$ be the $j$-th component of $L_I^\Phi(A,B)$.
Fix $(q_1,\dots,q_n) \in \cS^2(\RR^n)^n$
 and we will find $A,B \in M_n(\RR)$ such that
\begin{align}
\label{eqn:HAB 1}
q^{A,B}_j(e_j,e_j) & = q_j(e_j,e_j)\\
\label{eqn:HAB 2}
\Pair{e_i}{q^{A,B}_j(e_i,e_i)} + \Pair{e_j}{q^{A,B}_i(e_j,e_j)} 
 & = \Pair{e_i}{q_j(e_i,e_i)} + \Pair{e_j}{q_i(e_j,e_j)} \\
\label{eqn:HAB 3}
\Pair{e_1}{q_1^{A,B}(e_j,e_j)}
 & = \Pair{e_1}{q_1(e_j,e_j)}.
\end{align}
These equations imply 
 that $(q_1,\dots,q_n) -L_I^\Phi(A,B)$ is an element of $W$.

Take $A=(a_{ij}), B=(b_{ij}) \in M_n(\RR)$.
A direct computation with Lemma \ref{lemma:inj 0} implies that
\begin{align}
q^{A,B}_j(e_j,e_j)
 & = A \circ Q_{e_j}(e_j,e_j)-2Q_{e_j}(Ae_j,e_j)+Q_{Be_j}(e_j,e_j) \nonumber\\
 & = Ae_j+\sum_{k=1}^nb_{kj}Q_{e_k}(e_j,e_j)\nonumber\\
\label{eqn:HAB 4}
 & =(a_{jj}-b_{jj})e_j+\sum_{k \neq j}(a_{kj}+b_{kj})e_k,\\
q^{A,B}_i(e_j,e_j)
 & = A \circ Q_{e_i}(e_j,e_j)-2Q_{e_i}(Ae_j,e_j)+Q_{Be_i}(e_j,e_j)\nonumber\\
 & = Ae_i-2\sum_{k=1}^na_{kj}Q_{e_i}(e_k,e_j)
  +\sum_{k=1}^nb_{ki}Q_{e_k}(e_j,e_j)\nonumber\\
 & =(a_{ii}-2a_{jj}+b_{ii})e_i
 +(a_{ji}+2a_{ij}-b_{ji})e_j
 +\sum_{k \neq i,j}(a_{ki}+b_{ki})e_k.\nonumber
\end{align}
 for any mutually distinct $i,j=1,\dots,n$.
The latter equation implies that
\begin{equation}
\label{eqn:HAB 5}
\Pair{e_i}{q^{A,B}_j(e_i,e_i)} +\Pair{e_j}{q^{A,B}_i(e_j,e_j)}
  = 3(a_{ij}+a_{ji})-(b_{ij}+b_{ji})
\end{equation}
 for any mutually distinct $i,j=1,\dots,n$ and
\begin{equation}
\label{eqn:HAB 6}
 \Pair{e_1}{q^{A,B}_1(e_j,e_j)} = a_{11}-2a_{jj}+b_{11}
\end{equation}
 for any $j=2,\dots,n$.
 
Put $s_{ij}=\Pair{e_i}{q_j(e_j,e_j)}$,
 $t_{ij}=\Pair{e_j}{q_i(e_j,e_j)}$,
 and $u_j=\Pair{e_1}{q_1(e_j,e_j)}$
 for $i,j=1,\dots,n$.
Remark that $s_{11}=t_{11}=u_1$.
Put
 $a_{11}=s_{11}/2$, $b_{11}=-s_{11}/2$,
\begin{equation*}
a_{jj}=-u_j/2,\hsp b_{jj}=-s_{jj}-(u_j/2)
\end{equation*}
 for $j=2,\dots,n$, and
\begin{align*}
a_{ij} & = \frac{1}{4}(s_{ij}+t_{ij}),\\
b_{ij} & = s_{ij} - a_{ij} = \frac{1}{4}(3s_{ij}-t_{ij})
\end{align*}
 for any mutually distinct $i,j=1,\dots,n$.
By the equations
 (\ref{eqn:HAB 4}), (\ref{eqn:HAB 5}), and (\ref{eqn:HAB 6}),
 $A=(a_{ij})$ and $B=(b_{ij})$ satisfy
 the equations
 (\ref{eqn:HAB 1}), (\ref{eqn:HAB 2}), and (\ref{eqn:HAB 3}).
\end{proof}

Fix $(q_1,\dots,q_n) \in \Ker L_I^\Psi \cap  W$.
By the lemma, the goal is to show that $q_1=\cdots=q_n=0$.
\begin{lemma}
\label{lemma:inj 1}
 $q_j(e_i,e_j)=q_j(e_j,e_i)=0$ for any $i,j=1,\dots,n$.
\end{lemma}
\begin{proof}
When $i=j$, it is just shown by Equation (\ref{eqn:H 1})
 in the definition of $W$.
Take mutually distinct $i,j=1,\dots,n$.
Then,
\begin{align*}
0 & = \frac{1}{3}([q_i,Q_{e_j}]-[q_j,Q_{e_i}])(e_j,e_j,e_j) \\
 & =\left\{q_i(e_j,Q_{e_j}(e_j,e_j))-Q_{e_j}(e_j,q_i(e_j,e_j))\right\}\\
 & \quad -\left\{q_j(e_j,Q_{e_i}(e_j,e_j))-Q_{e_i}(e_j,q_j(e_j,e_j))\right\}\\
 & = \left\{q_i(e_j,-e_j)+q_i(e_j,e_j)\right\}
 -\left\{q_j(e_j,e_i)-Q_{e_i}(e_j,0)\right\}\\
 & = -q_j(e_j,e_i).
\end{align*}
Since $q_j$ is symmetric, we also obtain that $q_j(e_i,e_j)=0$.
\end{proof}
\begin{lemma}
\label{lemma:inj 2}
For any $i,j=1,\dots,n$,
\begin{equation}
\label{eqn:inj 2-1}
\Pair{e_i}{q_i(e_j,e_j)}+\Pair{e_j}{q_j(e_i,e_i)}=0.
\end{equation}
For any $i,j,k=1,\dots,n$ with $i \neq k$,
\begin{equation}
\label{eqn:inj 2-2}
\Pair{e_k}{q_i(e_j,e_j)}=0.
\end{equation}
\end{lemma}
\begin{proof}
When $i=j$, Lemma follows from the definition of $W$.
Suppose that $i \neq j$.
Since $q_i(e_i,e_j)=q_j(e_i,e_j)=q_j(e_j,e_j)=0$
 by Lemma \ref{lemma:inj 1} and
 Equation (\ref{eqn:H 1}) in the definition of $W$,
 we have
\begin{align*}
 [q_i,Q_{e_j}](e_i,e_j,e_j)
 & = \left\{q_i(e_i,Q_{e_j}(e_j,e_j))+2q_i(e_j,Q_{e_j}(e_i,e_j))\right\}\\
 & \quad - \left\{Q_{e_j}(e_i,q_i(e_j,e_j))+2Q_{e_j}(e_j,q_i(e_i,e_j))\right\}\\
 & = \left\{q_i(e_i,-e_j)+2q_i(e_j,-e_i)\right\}\\
 & \quad -\left\{\Pair{e_i}{q_i(e_j,e_j)} \cdot e_j
  -\Pair{e_j}{q_i(e_j,e_j)} \cdot e_i+2Q_{e_j}(e_j,0)\right\}\\
 & = \Pair{e_j}{q_i(e_j,e_j)} \cdot e_i-\Pair{e_i}{q_i(e_j,e_j)} \cdot e_j,\\
 [q_j,Q_{e_i}](e_i,e_j,e_j)
 &= \left\{q_j(e_i,Q_{e_i}(e_j,e_j))+2q_j(e_j,Q_{e_i}(e_i,e_j))\right\}\\
 & \quad -\left\{Q_{e_i}(e_i,q_j(e_j,e_j))+2 Q_{e_i}(e_j,q_j(e_i,e_j))\right\}\\
 &= \left\{q_j(e_i,e_i)+2q_j(e_j,-e_j)\right\}
  -\left\{Q_{e_i}(e_i,0)+2Q_{e_i}(e_j,0)\right\}\\
 & =q_j(e_i,e_i).
\end{align*}
Since $[q_i,Q_{e_j}]-[q_j,Q_{e_j}]=0$,
\begin{equation*}
q_j(e_i,e_i) =
\Pair{e_j}{q_i(e_j,e_j)} \cdot e_i-\Pair{e_i}{q_i(e_j,e_j)} \cdot e_j.
\end{equation*}
By taking the inner product with $e_k$, 
 we obtain that
 $\Pair{q_i(e_j,e_j)}{e_k}=0$ for $k \neq i,j$.
By taking the inner product with $e_i$ and $e_j$, 
 we also have
\begin{gather*}
 \Pair{e_i}{q_j(e_i,e_i)}-\Pair{e_j}{q_i(e_j,e_j)}=0\\
 \Pair{e_j}{q_j(e_i,e_i)}+\Pair{e_i}{q_i(e_j,e_j)}=0.
\end{gather*}
The latter is Equation (\ref{eqn:inj 2-1}).
Equation (\ref{eqn:inj 2-2}) follows from
 the former and Equation (\ref{eqn:H 2}) in the definition of $W$.
\end{proof}

Equations (\ref{eqn:H 3}) and (\ref{eqn:inj 2-1}) imply that
\begin{equation}
\label{eqn:inj 2-3}
 \Pair{e_1}{q_1(e_j,e_j)}=\Pair{e_j}{q_j(e_1,e_1)}=0.
\end{equation}
 for any $j=1,\dots,n$.
Now, we prove Proposition \ref{prop:vanishing} for $n=2$.
\begin{prop}
If $n=2$, then $\Ker L_I^\Psi=\Img L_I^\Phi$. 
\end{prop}
\begin{proof}
For $(q_1,q_2) \in \Ker L_I^\Psi \cap W$,
 $\Pair{e_i}{q_j(e_k,e_l)}=0$ for any $i,j,k,l=1,2$
 by Lemmas \ref{lemma:inj 1}, \ref{lemma:inj 2}
 and Equation (\ref{eqn:inj 2-3}).
Therefore, $q_1=q_2=0$.
Lemma \ref{lemma:H} implies that $\Ker L_I^\Psi=\Img L_I^\Phi$.
Proposition \ref{prop:vanishing} for $n=2$
 follows from Lemma \ref{lemma:standard}.
\end{proof}

We continue the proof for $n \geq 3$.
\begin{lemma}
\label{lemma:inj 3} 
$q_i(e_j,e_k)=q_j(e_k,e_i)=q_k(e_i,e_j)$
 for mutually distinct $i,j,k=1,\dots,n$.
\end{lemma}
\begin{proof}
Since $i,j,k$ are mutually distinct,
 Lemma \ref{lemma:inj 2} implies
\begin{align*}
\frac{1}{3} \cdot [q_i,Q_{e_j}](e_k,e_k,e_k)
 & = q_i(e_k,Q_{e_j}(e_k,e_k))-Q_{e_j}(e_k,q_i(e_k,e_k))\\
 & =q_i(e_k,e_j)
  -\left\{\Pair{e_k}{q_i(e_k,e_k)}\cdot e_j
  -\Pair{e_j}{q_i(e_k,e_k)} \cdot e_k \right\}\\
 & = q_i(e_k,e_j).
\end{align*}
Similarly, we have $(1/3)\cdot [q_j,Q_{e_i}](e_k,e_k,e_k)=q_j(e_k,e_i)$.
Hence,
\begin{equation*}
 q_i(e_k,e_j)-q_j(e_k,e_i)
 = \frac{1}{3} \cdot ([q_i,Q_{e_j}]-[q_j,Q_{e_i}])(e_k,e_k,e_k)=0.
\end{equation*}
It implies $q_i(e_j,e_k)=q_i(e_k,e_j)=q_j(e_k,e_i)$.
By permutations of indices $(i,j,k)$,
 we obtain that $q_j(e_k,e_i)=q_k(e_i,e_j)$.
\end{proof}

\begin{lemma}
\label{lemma:inj 4} 
For $i,j,k=1,\dots,n$,
\begin{gather}
\label{eqn:inj 4-1}
 q_i(e_j,e_j)=0,\\
\label{eqn:inj 4-2}
\Pair{e_i}{q_i(e_j,e_k)}=\Pair{e_j}{q_i(e_j,e_k)}=\Pair{e_k}{q_i(e_j,e_k)}=0.
\end{gather}
\end{lemma}
\begin{proof}
For mutually distinct $i,j,k=1,\dots,n$,
\begin{align*}
[q_i,Q_{e_j}](e_j,e_k,e_k)
 & = \left\{q_i(e_j,Q_{e_j}(e_k,e_k))+2q_i(e_k,Q_{e_j}(e_j,e_k))\right\}\\
 & \quad - \left\{Q_{e_j}(e_j,q_i(e_k,e_k))+2Q_{e_j}(e_k,q_i(e_j,e_k))\right\}\\
 & = \left\{q_i(e_j,e_j)+2q_i(e_k,-e_k)\right\}\\
 & \quad -\left\{
 -q_i(e_k,e_k)+2\left(\Pair{e_k}{q_i(e_j,e_k)} \cdot e_j
 -\Pair{e_j}{q_i(e_j,e_k)} \cdot e_k\right)\right\}\\
 & = q_i(e_j,e_j)-q_i(e_k,e_k)\\
 & \quad -2\Pair{e_k}{q_i(e_j,e_k)} \cdot e_j
 +2\Pair{e_j}{q_i(e_j,e_k)} \cdot e_k, \\
[q_j,Q_{e_i}](e_j,e_k,e_k)
 & = \left\{q_j(e_j,Q_{e_i}(e_k,e_k))+2q_j(e_k,Q_{e_i}(e_j,e_k))\right\}\\
 & \quad - \left\{Q_{e_i}(e_j,q_j(e_k,e_k))+2Q_{e_i}(e_k,q_j(e_j,e_k))\right\}\\
 & = \left\{q_j(e_j,e_i)+2q_j(e_k,0)\right\}\\
 &\quad -\left\{\Pair{e_j}{q_j(e_k,e_k)}\cdot e_i
 -\Pair{e_i}{q_j(e_k,e_k)}\cdot e_j+2Q_{e_i}(e_k,0) \right\}\\
 & = -\Pair{e_j}{q_j(e_k,e_k)}\cdot e_i.
\end{align*} 
Since $[q_i,Q_{e_j}]-[q_j,Q_{e_i}]=0$, we obtain that
\begin{equation*}
 q_i(e_j,e_j)-q_i(e_k,e_k)
 =-\Pair{q_j(e_k,e_k)}{e_j}\cdot e_i+2\Pair{q_i(e_j,e_k)}{e_k} \cdot e_j
 -2\Pair{q_i(e_j,e_k)}{e_j} \cdot e_k.
\end{equation*}
By taking the inner product of the with $e_i$ and $e_j$,
\begin{align}
\Pair{e_i}{q_i(e_j,e_j)}-\Pair{e_i}{q_i(e_k,e_k)}
 & = - \Pair{e_j}{q_j(e_k,e_k)},\nonumber\\
\label{eqn:inj 4-3}
\Pair{e_j}{q_i(e_j,e_j)}-\Pair{e_j}{q_i(e_k,e_k)}
 & = 2\Pair{e_k}{q_i(e_j,e_k)}.
\end{align}
The former equation for $i=1$
 implies $\Pair{e_j}{q_j(e_k,e_k)}=0$
 for any mutually distinct $j,k=2,\dots,n$.
By Equation (\ref{eqn:inj 2-3}),
 the same equation holds for the case $j=1$ or $k=1$.
Combined with Equation (\ref{eqn:inj 2-2}),
 we obtain Equation (\ref{eqn:inj 4-1}).

Equations (\ref{eqn:inj 2-2}) and (\ref{eqn:inj 4-3}) imply
 $\Pair{e_k}{q_i(e_j,e_k)}=0$.
By permutations of indices $(i,j,k)$ and Lemma \ref{lemma:inj 3},
 we obtain Equation (\ref{eqn:inj 4-2})
 for mutually distinct $i,j,k$.
Equation (\ref{eqn:inj 4-2}) for other cases
 follows from Lemma \ref{lemma:inj 1}
 and Equation (\ref{eqn:inj 4-1}).
\end{proof}

Proposition \ref{prop:vanishing} for $n=3$ follows from the lemma.
\begin{prop}
If $n=3$, then $\Ker L_I^\Psi= \Img L_I^\Phi$.
\end{prop}
\begin{proof}
For $(q_1,q_2,q_3) \in \Ker L_I^\Psi \cap W$,
 Equation (\ref{eqn:inj 4-2}) in Lemma \ref{lemma:inj 4} implies
 $q_1=q_2=q_3=0$ if $n=3$.
By Lemma \ref{lemma:H}, we have $\Ker L_I^\Psi= \Img L_I^\Phi$.
Proposition \ref{prop:vanishing} for $n=3$
 follows from Lemma \ref{lemma:standard}.
\end{proof}

The following lemma completes the proof for $n \geq 4$.
\begin{lemma}
\label{lemma:inj 5} 
$q_i(e_j,e_k)=0$ for any $i,j,k=1,\dots,n$.
\end{lemma}
\begin{proof}
By Lemma \ref{lemma:inj 1} and \ref{lemma:inj 4},
 it is sufficient to show that $\Pair{e_i}{q_j(e_k,e_l)}=0$
 for mutually distinct $i,j,k,l=1,\dots,n$.
Take mutually disjoint $i,j,k,l=1,\dots,n$.
Then,
\begin{align*}
[q_i,Q_{e_j}](e_k,e_l,e_l)
 & = \left\{q_i(e_k,Q_{e_j}(e_l,e_l))+2q_i(e_l,Q_{e_j}(e_k,e_l))\right\}\\
 & \quad -\left\{Q_{e_j}(e_k,q_i(e_l,e_l))+2Q_{e_j}(e_l,q_i(e_k,e_l))\right\}\\
 & = \left\{q_i(e_k,e_j)+2q_i(e_l,0)\right\}\\
 & \quad -\left\{Q_{e_j}(e_k,0)
 + 2 (\Pair{e_l}{q_i(e_k,e_l)} \cdot e_j
 -\Pair{e_j}{q_i(e_k,e_l)} \cdot e_l)\right\}\\
 & = q_i(e_j,e_k)-2 \Pair{e_j}{q_i(e_k,e_l)} \cdot e_l.
\end{align*}
Similarly, we obtain that
\begin{equation*}
 [q_j,Q_{e_i}](e_k,e_l,e_l)
 = q_j(e_i,e_k)-2 \Pair{e_i}{q_j(e_k,e_l)} \cdot e_l.
\end{equation*}
Since $[q_i,Q_{e_j}]-[q_j,Q_{e_i}]=0$,
\begin{equation*}
 q_i(e_j,e_k)-q_j(e_i,e_k)
 = \left\{\Pair{e_j}{q_i(e_k,e_l)}-\Pair{e_i}{q_j(e_k,e_l)}\right\} \cdot e_l.
\end{equation*}
By Lemma \ref{lemma:inj 3}, $q_i(e_j,e_k)=q_j(e_k,e_l)$
 and $q_i(e_k,e_l)=q_k(e_l,e_i)$.
Hence, we have
\begin{equation*}
 \Pair{e_j}{q_k(e_l,e_i)}=\Pair{e_i}{q_j(e_k,e_l)}.
\end{equation*}
By take permutations of indices $(i,j,k,l)$,
\begin{equation}
\label{eqn:inj 5-1}
 \Pair{e_l}{q_i(e_j,e_k)}=\Pair{e_i}{q_j(e_k,e_l)}
  =\Pair{e_j}{q_k(e_l,e_i)}=\Pair{q_l(e_i,e_j)}{e_k}.
\end{equation}
On the other hand, we have
\begin{align*}
 [q_i,Q_{e_j}](e_j,e_k,e_l)
 & =\left\{q_i(e_j,Q_{e_j}(e_k,e_l))+q_i(e_k,Q_{e_j}(e_l,e_j))
 +q_i(e_l,Q_{e_j}(e_j,e_k)) \right\} \\
 & \quad -\left\{Q_{e_j}(e_j,q_i(e_k,e_l))+Q_{e_j}(e_k,q_i(e_l,e_j))
 +Q_{e_j}(e_l,q_i(e_j,e_k))\right\}\\
 & = \left\{q_i(e_j,0)+q_i(e_k,-e_l)+q_i(e_l,-e_k)\right\}\\
 & \quad -\left\{
 -q_i(e_k,e_l)+\Pair{e_k}{q_i(e_l,e_j)} \cdot e_j
 +\Pair{e_l}{q_i(e_j,e_k)} \cdot e_j\right\}\\
 & =-q_i(e_k,e_l)-2 \Pair{e_j}{q_i(e_k,e_l)} \cdot e_j,
\end{align*}
 and
\begin{align*}
[q_j,Q_{e_i}](e_j,e_k,e_l)
 & =\left\{q_j(e_j,Q_{e_i}(e_k,e_l))+q_j(e_k,Q_{e_i}(e_l,e_j))
 +q_j(e_l,Q_{e_i}(e_j,e_k)) \right\} \\
 & \quad -\left\{Q_{e_i}(e_j,q_j(e_k,e_l))+Q_{e_i}(e_k,q_j(e_l,e_j))
 +Q_{e_i}(e_l,q_j(e_j,e_k))\right\}\\
 & =\left\{q_j(e_j,0)+q_j(e_k,0)+q_j(e_l,0)\right\}\\
 & \quad -\left\{-\Pair{e_i}{q_j(e_k,e_l)} \cdot e_j
 +Q_{e_i}(e_k,0)+Q_{e_i}(e_l,0) \right\}\\
 & =\Pair{e_j}{q_i(e_k,e_l)} \cdot e_j.
\end{align*}
Since $[q_i,Q_{e_j}]-[q_j,Q_{e_i}]=0$, 
\begin{equation*}
 q_i(e_k,e_l)+3 \cdot \Pair{e_j}{q_i(e_k,e_l)}\cdot e_j=0.
\end{equation*}
By taking the inner product with $e_j$,
 we have $\Pair{e_j}{q_i(e_k,e_l)}=0$.
Hence,
\begin{equation*}
\Pair{e_i}{q_j(e_k,e_l)}=0 
\end{equation*}
 by permuting indices $(i,j,k,l)$.
\end{proof}

Now, we prove Proposition \ref{prop:vanishing} for $n \geq 4$.
The last lemma implies that $q_1=\dots=q_n=0$
 for any $(q_1,\dots,q_n) \in \Ker L_I^\Psi \cap W$.
By Lemma \ref{lemma:H}, we obtain that $\Ker L_I^\Psi=\Img L_I^\Phi$.
Proposition \ref{prop:vanishing}
 follows from Lemma \ref{lemma:standard}.
\bigskip

\section{Classification of the standard actions}
In this section, we classify the standard $\Gamma_{n,k}$-actions up
 to smooth conjugacy.
Let $O(n)$ be the orthogonal group of $\RR^n$.
\begin{prop}
\label{prop:classify} 
For $B,B' \in \GL_n(\RR)$,
 $\rho_B$ and $\rho_{B'}$ are smoothly conjugate
 if and only if 
 there exists $T \in O(n)$ and $c>0$ such that $B'= (cT)B$.
\end{prop}
Remark that all standard $\Gamma_{n,k}$-actions are 
 {\it topologically conjugate} to each other
 {\it i.e.} there exists a {\it homeomorphism} $h$ of $S^n$
 such that $\rho_{B'}^\gamma \circ h=h \circ \rho_B^\gamma$
 for any $\gamma \in \Gamma_{n,k}$.
In fact, if $B'=AB$ for some $A \in \GL_n(\RR)$,
 then the linear map $x \mapsto Ax$ on $\RR^n$ extends
 to a {\it homeomorphism} $h_A$ on $S^n$.
It is easy to check that
 $\rho_{B'}^\gamma \circ h_A=h_A \circ \rho_B^\gamma$
 for any $\gamma = a,b_1,\dots,b_n$.
When $A=cT$ with $c>0$ and $T \in O(n)$,
 then $h_A$ is a {\it diffeomorphism}.
Hence, $\rho_B$ and $\rho_{B'}$ are smoothly conjugate in this case.

To prove the ``only if'' part of Proposition \ref{prop:classify},
 we need a technical lemma.
Recall that $Q_v \in \cS^2(\RR^n)$ is defined by
\begin{equation}
Q_v(\xi,\eta)
 =\Pair{\xi}{\eta} \cdot v-\Pair{\xi}{v} \cdot \eta
 -\Pair{\eta}{v} \cdot \xi.
\end{equation}
\begin{lemma}
\label{lemma:commutes} 
Suppose that $A \circ Q_v=Q_w \circ (A,A)$
 for $v,w \in \RR^n \sm \{0\}$ and $A \in \GL_n(\RR)$.
Then, $A=cT$ for some $c>0$ and $T \in O(n)$.
\end{lemma}
\begin{proof}
By a direct computation, we have
\begin{align*}
Q_w(A\xi,A\xi) &=\|A\xi\|^2 \cdot w -2 \Pair{A\xi}{w} \cdot A\xi,\\
A \circ Q_v(\xi,\xi) & =\|\xi\|^2 \cdot Av -2 \Pair{\xi}{v} \cdot A \xi.
\end{align*}
Hence,
\begin{equation}
\label{eqn:commutes 1}
 \|A\xi\|^2 \cdot w - \|\xi\|^2 \cdot Av
 =2\left(\Pair{A\xi}{w}-\Pair{\xi}{v}\right) \cdot A\xi
\end{equation}
 for any $\xi \in \RR^n$.
Put $\lambda=\left(2\Pair{Av}{w}-\|v\|^2 \right)/\|Av\|^2$.
Then, the equation for $\xi=v$ implies $w=\lambda Av$.
By substituting it to Equation (\ref{eqn:commutes 1}), we have
\begin{equation*}
(\lambda \|A\xi\|^2-\|\xi\|^2)Av
 =2\left(\lambda \Pair{A\xi}{Av}-\Pair{\xi}{v}\right) \cdot A\xi.
\end{equation*}
Since $A$ is invertible,
 it implies that $\|A\xi\|=\lambda^{-1}\|\xi\|$
 for any $\xi \in \RR^n \sm \RR v$.
Since $\RR^n \sm \RR v$ is a dense subset of $\RR^n$,
 the same holds for any $\xi \in \RR^n$.
Hence, there exists $T \in O(n)$ such that $A=\lambda^{-1}T$.
\end{proof}

\begin{proof}
[Proof of Proposition \ref{prop:classify}]
It is sufficient to show the ``only if'' part.
Suppose that $\rho_B$ and $\rho_{B'}$ are smoothly conjugate.
Take a diffeomorphism $h$ of $S^n$ such that
 $\rho_{B'}^\gamma \circ h =h \circ \rho_B^\gamma$
 for any $\gamma \in \Gamma_{n,k}$.
Since $\infty$ is the unique global fixed point of $\rho_B$ and $\rho_{B'}$,
 the diffeomorphism $h$ fixes $\infty$.
Recall that $P_B$ and $P_{B'}$ are the local $\Gamma_{n,k}$-actions
 defined by
 $P_B^\gamma=\bar{\phi} \circ \rho_B^\gamma \circ \bar{\phi}^{-1}$
 and $P_{B'}^\gamma=\bar{\phi} \circ \rho_{B'}^\gamma \circ \bar{\phi}^{-1}$,
 where $\bar{\phi}(x)=(1/\|x\|^2) \cdot x$.
Put $H=\bar{\phi} \circ h \circ \bar{\phi}^{-1}$
 and $A=D^{(1)}_0 H$.
Then, $P_{B'}^\gamma \circ H=H \circ P_B^\gamma$,
 and hence,
\begin{equation}
\label{eqn:classify 1}
 D^{(1)}_0 P_{B'}^\gamma \circ D^{(2)}_0 H
 +D^{(2)}_0 P_{B'}^\gamma \circ (A,A)
 = A \circ D^{(2)}_0 P_B^\gamma
 + D^{(2)}_0H \circ (D^{(1)}_0P_B^\gamma,D^{(1)}_0 P_B^\gamma).
\end{equation}
Since $P_B^a(x) = P_{B'}^a(x)=k^{-1}x$,
  the equation for $\gamma=a$ implies $k^{-1} D^{(2)}_0 H= k^{-2}D^{(2)}_0H$.
Therefore, $D^{(2)}_0 H=0$.
Put $B=(v_1,\dots,v_n)$ and $B'=(w_1,\dots,w_n)$.
Since $D^{(2)}_0 P_B^{b_i}=2Q_{v_i}$ and $D^{(2)}_0 P_{B'}^{b_i}=2Q_{w_i}$,
Equation (\ref{eqn:classify 1}) for $\gamma=b_i$ implies
\begin{equation*}
 Q_{w_i} \circ (A,A)=A \circ Q_{v_i}
\end{equation*}
 for any $i=1,\dots,n$.
By Lemma \ref{lemma:commutes},
 there exists $c>0$ and $T \in O(n)$ such that $A=cT$.
Since $T$ preserves the inner product,
\begin{align*}
(cT) \circ Q_{v_i}(\xi,\eta)
 & =Q_{w_i}((cT) \xi,(cT) \eta) \\
 & = c^2 \left\{\Pair{T \xi}{T \xi}\cdot w_i -\Pair{T \xi}{w_i} \cdot T\eta
   -\Pair{T \eta}{w_i} \cdot T \xi
 \right\}\\
 & =(cT) \circ \left\{\Pair{\xi}{\eta} \cdot (cT^{-1})w_i
  -\Pair{\xi}{(cT^{-1})w_i} \cdot \eta
  -\Pair{\eta}{(cT^{-1})w_i} \cdot \xi \right\}\\
 & = (cT) \circ Q_{cT^{-1}w_i}(\xi,\eta)
\end{align*}
 for any $\xi,\eta \in \RR^n$.
It implies that $v_i=cT^{-1}w_i$ for any $i=1,\dots,n$.
Therefore, $B'=(c^{-1}T)B$.
\end{proof}

\end{document}